\hsize=30pc
\vsize=50pc

\ifx\shlhetal\undefinedcontrolsequence\let\shlhetal\relax\fi
\input amstex.tex
\documentstyle{amsppt}
  \catcode`@11

\catcode`@12

\def\martinstricks{}   

\newtoks\fooline

\fooline={}

\catcode`@11

\newtoks\oldoutput
\oldoutput=\the\output
\output={\setbox255\vbox to \ht255 
         {\unvbox255 %
\iffirstpage@\else 
\smallskip\line{\the\fooline}\vss\fi}\the\oldoutput}

\catcode`@12

\martinstricks

\def\demo#1{\smallskip\noindent{\smc #1}}
\def\enddemo{\medskip}

\def\lec{\mathrel{\mathord{<}\!\!\raise 0.8 pt\hbox{$\scriptstyle\circ$}}}

\newcount\skewfactor
\def\mathunderaccent#1#2 {\let\theaccent#1\skewfactor#2
\mathpalette\putaccentunder}
\def\putaccentunder#1#2{\oalign{$#1#2$\crcr\hidewidth
\vbox to.2ex{\hbox{$#1\skew\skewfactor\theaccent{}$}\vss}\hidewidth}}
\def\name{\mathunderaccent\tilde-3 }

\NoBlackBoxes
\documentstyle {amsppt}
\topmatter
\title {{\it On some problems in general topology} \\
Sh-E3} \endtitle
\subjclass{54A25, 54D15}\endsubjclass

\author {Saharon Shelah \thanks {
Partially supported by the NSF. 
Research in this paper was done
in fall 77 (while the author was at University of Wisconsin, Madison),
typed and distributed in winter 78 (while the author was at the
University of California, Berkeley), but the author
forgot to submit it for publication.}
\endthanks} \endauthor

\endtopmatter
\document


\head{ \S 0. Introduction }\endhead

This work was done in 1977 and was widely quoted but not submited.

In section  3 it is proved that Arhangelskii's problem has a consistent
positive answer: if $V\vDash CH$, then for some $\aleph_1$-complete
$\aleph_2$-c.c. forcing notion $P$ of cardinality $\aleph_2$ we have 
$\Vdash_P$ ``$CH$ and there is a
Lindel\"of regular topological space of size $\aleph_2$ with clopen
basis with every point of 
pseudo-character $\aleph_0$ (i.e. each singleton is the intersection of
countably many open sets)''.

Meanwhile this was continued in Hajnal and Juhasz [HJ],
Stanley and Shelah [ShSt:167], I.Gorelic [Go] and Ch.Morgan.

In section 4 we prove the consistency of: $CH+ 2^{\aleph_1} > \aleph_2 +$
there is no space as above with $\aleph_2$ points" (starting with a
weakly compact cardinal).

Section 2 deals with $\beta(\Bbb N)$, it is proved that the
following is consistent with $ZFC$: $MA+2^{\aleph_0} = \aleph_2 + (*)$
where

\line{$(*)$\hfil \vtop{\hsize=0.95\hsize \noindent
 if $A^0_i,A^1_i\subseteq \omega\,
(\text{ for }  i<\omega_1)$ and $A^0_i\cap A^1_j$ is finite
for $i,j<\omega_1$; and $\goth D^\ell_i$ is a non-principal
ultrafilter over $\omega$ such that $A^\ell_i\in \goth
D^\ell_i\,(\text{ for }i<\omega_1\text{ and }\ell\in
\{0,1\})$, {\it then} there is a
$B\subseteq \omega$ such that $B\in \goth D^\ell_i
\Leftrightarrow \ell=0$.}}
  \bigskip

The scheme is as in Baumgartner [B].

Another problem on $\beta(\Bbb N)$ which I remember was asked and
published by E. van Douwen and G. Woods, is answered in \S1: is there
a discrete $D\subseteq \beta(\Bbb N)$, of cardinality $\aleph_1$ and
$A\subseteq D$ such that $cl(A) \cap cl(D\setminus A)\neq \emptyset$.

I thank U. Abraham for urging the publication (in this form) and for
corrections,  the referee for
corrections 
and M. D\v{z}amonja for corrections.
Compared with the old version we added details, explanations, the
introduction, added 2.4 and stated 1.2, 4.2;
one section was omitted.

\head{\S 1. A problem on $\beta(\Bbb N)$}\endhead

\subhead{1.1 Question}\endsubhead  Does there exist a discrete 
$D \subseteq\beta(\Bbb N)$ of cardinality $\aleph_1$, and
an $A\subseteq  D$ such that $cl(A) \cap
cl( D\setminus A)\not=\emptyset?$ 
\bigskip
\subhead{1.2 Answer}\endsubhead  Yes, moreover we can let $ D$
have any cardinality  $\lambda$ such that $\aleph_1\leq
\lambda= cf(\lambda)\le 2^{\aleph_0}$ and have $\bigcap\{cl(D'):
D'\subseteq D\text{ and }|D'|=|D|\}$ not empty.
\bigskip
\subhead{1.3 Definition}\endsubhead  Let $B_i$ (for $i\leq\omega_1$) be
the $BA$ (Boolean algebra) freely generated by
$\{x_\alpha:\alpha <i\}, B^c_i$ is the completion of $B_i$,
$B=B^c_{\omega_1}$. 
\bigskip

\proclaim{1.4 Claim}
 In the space of ultrafilters of
$B^c_{\omega_1}$ we can find such a $D$. \endproclaim
\bigskip
\demo{PROOF:}  We define by induction on $i\leq \omega_1$,
an ultrafilter $\goth D_i$ of $B^c_{\omega_1}$ such that
\roster
\item"{(i)}" $x_\alpha \in
\goth D_i\Leftrightarrow\alpha=i$, 
\item"{(ii)}"  if $a\in
B^c_i$, $a\in \goth D_i$ then $a\in \goth D_j$ for every
$j>i$.  
\endroster
This is easy. Let $\goth D_{\omega_1}= \bigcup_{i<\omega_1} (\goth D_i
\cap B^c_i)$.
Now, $\{\goth D_i:i<\omega_1\}$ is discrete
by (i) and $\goth D_{\omega_1} \in cl\{\goth D_i:i\in
S\}$ for any $S\subseteq \omega_1, |S|=\aleph_1$, because
$B^c_{\omega_1}=\bigcup_{i<{\omega_1}} B^c_i$, as
$B_{\omega_1}$ satisfies the countable chain condition.
\enddemo
\bigskip
\subhead{1.5 Solution of the problem}\endsubhead  Let
$X_i\subseteq \omega \  (i<\omega_1)$ be independent, i.e.
any non trivial Boolean combination of the $X_i$ is not
empty.  Let $f:\Cal P(\omega)\rightarrow B^c_{\omega_1}$ be
any homomorphism such that $f(X_i)=x_i$ (exists as
$B^c_{\omega_1}$ is complete and the $X_i$ are
independent).  It is not hard to prove $f$ is onto, and
$\{f^{-1}(\goth D_i): i<\omega_1\}$ is as required $(\goth
D_i$ from the claim), and
$$f^{-1}(\goth D_{\omega_1})\in cl(\{f^{-1}(\goth
D_{2i}):i<\omega_1\}) \cap cl(\{f^{-1}(\goth D_{2i+1}):i<
 \omega_1\}).$$
\null\hfill$\square_{1.5}$
\bigskip

\head{\S 2 A question on $\beta(\Bbb N)\setminus \Bbb N$}\endhead

\proclaim{2.1 Claim} 
Assuming the consistency of $ZFC$ we prove the consistency of the
following assertion with 

$ZFC +2^{\aleph_0} = \aleph_2+MA$:

\indent $\otimes$  if $A^0_i,A^1_i\subseteq \omega\,
(\text{ for }  i<\omega_1)$ and $A^0_i\cap A^1_j$ is finite
for $i,j<\omega_1;\ \goth D^\ell_i$ is a non-principal
ultrafilter over $\omega$ such that $A^\ell_i\in \goth
D^\ell_i\,(\text{ for }i<\omega_1\text{ and }\ell\in
\{0,1\})$, {\it then} there is a
$B\subseteq \omega$ such that $B\in \goth D^\ell_i
\Leftrightarrow \ell=0.$ \endproclaim \bigskip  
\demo{PROOF:}
We assume $V\vDash 2^{\aleph_0}=\aleph_1 \ \wedge \
2^{\aleph_1}=\aleph_2$.
We repeat the proof of Solovay-Tennenbaum of $Con(ZFC +
2^{\aleph_0} =\aleph_2 + MA)$, similarly to Baumgartner
[B]. That is, we define by induction on $\alpha <\omega_2$ a
set of forcing conditions $P_\alpha$, increasing (under $\subseteq$
and even $\lec$) and
continuous with $\alpha, |P_\alpha|\leq \aleph_1$, and
each $P_\alpha$ satisfies the countable chain condition.  We
start with $V\vDash
2^{\aleph_0}=\aleph_1+2^{\aleph_1}=\aleph_2
+\diamondsuit_{\{\delta<{\omega_2}:cf(\delta)=\aleph_1\}}$.

 Now, in
addition, at some $\alpha<\omega_2$ we
consider a system $\langle A^\ell_i: i<\omega_1,
\ell=0,1\rangle$, $\langle \goth D^\ell_i : i
<\omega_1,\ell=0,1\rangle$, which belongs to $V^{P_\alpha},
A^\ell_i \subseteq \omega ,\  A^\ell_i\in \goth
D^\ell_i,\  \goth D^\ell_i$ a family of subsets of
$\omega$, such that the filter it generates in
$V^{P_\alpha}$ (which we denote by the same letter) is
$\aleph_1$-saturated (i.e. there are no $C_\alpha
\not=\emptyset\ mod\ \goth D^\ell_i,C_\alpha \subseteq
\omega$ for $\alpha<\omega_1$ such that $ C_\alpha \cap
C_\beta= \emptyset\ mod\ \goth D^\ell_i$ for $\alpha <\beta <
\omega_1)$.  By the usual bookkeeping, every such system
appears, and this is possible as $V\vDash 2^{\aleph_0}=\aleph_1 +
2^{\aleph_1}=\aleph_2$. Clearly $V^{P_\alpha}\vDash
2^{\aleph_0}=\aleph_1$.  We define in $V^{P_\alpha}$ a set
of forcing conditions $Q$ satisfying the c.c.c., 
whose generic set gives a $B$ such that $B\in \goth D^0_i$ and
$\omega\setminus B\in\goth D^1_i$  (for $i<\omega_1$), and define
$P_{\alpha+1}=P_\alpha * Q$.  This is sufficient, as if
$p\in P_{\omega_2}$ forces that
 $\langle \name A^\ell_i:
i<\omega_1\rangle, \langle\name{\goth D}^\ell_i: i
<\omega_1,\ell=0,1\rangle$ contradict $\otimes$, then for 
a club of $E$ of $\omega_2$, for every $\alpha\in E$ of
cofinality $\aleph_1, \langle
\name A^\ell_i:\ell <\omega_1\rangle$ is a $P_\alpha$-name,
and $\langle \Cal P(\omega)^{V^{P_\alpha}} \cap\name{\goth
D}^\ell_i:\ell, i\rangle$ is a $P_\alpha$-name of an ultrafilter. Then,
clearly for stationarily many $\alpha<\omega_2$ of
cofinality $\aleph_1$,
 in $V^{P_\alpha}$ the above holds and as
$\diamondsuit_{\{\delta<\omega_2: cf(\delta)=\aleph_1\}}$ holds we can
assume that
we have considered the system in question at some $\alpha$ (of course
in the bookkeeping we take care of $MA$ too).
\enddemo

So it suffices to prove:

\proclaim{2.2 Claim}  If $V\vDash 2^{\aleph_0}=\aleph_1$,
$A^\ell_i\subseteq
\omega\,(\text{ for } \ell=0,1$ and $i<\omega_1)$,
$A^0_i\cap A^1_j$ are finite for $i, j<\omega_1$, $\goth
D^\ell_i$ an $\aleph_1$-saturated filter over $\omega$,
$A^\ell_i\in \goth D^\ell_i$, {\it then} we can find a
partial order $Q$ of size $\aleph_1$, satisfying c.c.c., and in
$V^Q$ there is an $X\subseteq \omega$ such that $ X\in \goth
D^0_i$ and $\omega \setminus X\in \goth D^1_i$
(for $i<\omega_1$; note: $\goth D^\ell_i$ stands for the filter it
generates). 
\endproclaim \bigskip

\subhead{2.2A Remark}\endsubhead
The sequence $\langle (B^0_i, B^1_i): i<\omega_1\rangle$ constructed
below should be just generic enough (we do not use this presentation
because some people do not like it\footnote{at least they
did not back in the seventies}). More specifically
letting $f_{2i+\ell}\in {}^\omega 2$ be
defined by 
$$
f^\ell_i(n)=0 \Leftrightarrow [\text{the }n\text{th element of
}A^\ell_i\text{ belongs to }B^\ell_i]
$$
we demand: for every $n<\omega$ and open dense $J\subseteq {}^n(^\omega
2)$, for some $\alpha_J<\omega_1$, for every
$\alpha_0<\ldots<\alpha_{n-1}$ from $(\alpha_J, \omega_1)$ we have
$\langle f_{\alpha_0}, \ldots, f_{\alpha_{n-1}}\rangle \in J$.
\bigskip

\demo{PROOF:} We shall choose sets $B^\ell_i\subseteq
A^\ell_i$, $B^\ell_i\notin \goth D^\ell_i$, and let 
$Q=\{(f,g):f,g$ are finite functions from $\omega_1$ to
$\omega$, and 
  $(A^0_i\setminus B^0_i\setminus f(i))
\cap(A^1_j\setminus B^1_j\setminus g(j))=\emptyset$ when 
$i\in Dom (f)$ and $j\in Dom (g)\}$.

( let $Q_i= \{(f,g) \in Q: \ Dom (f)\cup Dom (g)\subseteq
i\}$).

$Q$ is ordered naturally: $(f_1, g_1)\leq (f_2, g_2)$ \underbar{iff} 
$f_1\subseteq f_2\ \&\ g_1\subseteq g_2$

For a generic $G\subseteq Q$ the set
$X=\bigcup\{B^0_i\setminus
f(i):\,\exists g\,((f,g)\in G)\}$ is as required,
and $|Q|=\aleph_1$.  We should show only that we can choose
$B^\ell_i$ such that $Q$ satisfies the c.c.c. (the density condition is easy: 
for every $i,j<\omega_1$ by the almost disjoint condition,
$A^0_i\cap A^1_j$ is finite, hence for every $(f,g)\in
Q$ and $i<\omega_1$ there is $n^*<\omega$ such that $A^0_i\subseteq
n^*$ is disjoint to $A^1_j$ for $j\in Dom(g)\cup \{i\}$ and
$A^1_i\setminus n^*$ is disjoint to $A^0_j$ for $j\in Dom(f)\cup
\{i\}$. Let $f'=f\cup \{\langle i, n^*\rangle\}$, $g'= g\cup \{\langle
i, n^*\rangle\}$.
So, there is a $(f',g')\in Q$ such
that $(f,g)\leq (f',g')$ and $i\in Dom (f')\cap Dom (g'))$.

Suppose $(f_i,g_i)\quad (\text{for } i<\omega_1)$ exemplify a
contradiction to c.c.c.  Then, by the well known techniques, we
can assume that there is a $(f,g) \leq (f_i,g_i)$,
$Dom(f_i)=Dom(g_i)=w_i$, $Dom (f)=Dom (g) =w$, $w\subseteq w_i$ and the
sets
$w_i\setminus w$ are pairwise disjoint,
$w_i\setminus w=\{\eta_i(0),\dots,
\eta_i(t)\}$, $i<j\Rightarrow\eta_i(t)< \eta_j(0)$,
$\eta_i(0)<\eta_i(1)<\ldots<\eta_i(t)$, and for $m\leq t$ we have
$f_i(\eta_i(m))=k_m$ and $g_i(\eta_i(m))=j_m$.  For some
$\delta$ we can get an ``elementary submodel'' of the whole
system (so $\eta_i(t)<\delta$ for $i<\delta)$.  Let
$t_\zeta,\delta_\zeta, f^\zeta_i,g^\zeta_i\, (\text{for }
i<\delta), \eta_i^\zeta(m), k^\zeta_m, j^\zeta_m
\,(\text{for} \  \zeta <\omega_1 \text{ and} \  m\leq t)$
enumerate all possible such systems.  Now we shall choose
$B^\ell_\alpha$ by induction on $\alpha$ and then on
$\ell$ with some restrictions:  (say $\ell=0$ for
notational simplicity).

Let us try to explain the idea of the proof.
If  
$\max\{\zeta, \delta_\zeta\}<\alpha<\omega_1$, then we think of $(\langle
\alpha, k^\zeta_0\rangle ,\emptyset)\in Q$  as a
candidate to be a  part of some $(f_i,g_i)$. Now, either it
is
compatible with infinitely many $(f^\zeta_i,g^\zeta_i)\,
(\text{for} \  i<\delta_\zeta)$, or the condition on being an 
elementary submodel eliminates  this possibility, and
similarly if 
$Dom (f')=Dom(g')=\{\eta(0),\dots,\eta(m-1)\}\subseteq \alpha$ and $(f,g)$ is
compatible with infinitely many $(f^\zeta_i,g^\zeta_i)$,
then either so is 
$(f'\cup \{\langle\alpha, k^\zeta_m\rangle\}, g')$, 
or the condition on elementary
submodels is violated. 

 There are countably many
such conditions and we can find $2^{\aleph_0}$
pairwise almost disjoint infinite $Y_\xi \subseteq A^\ell_i$   
(for $\xi<2^{\aleph_0})$, such that
$A^\ell_\alpha\setminus
 Y_\xi$ are as required and all
but countably many of the $Y_\xi$'s are
$=\emptyset \ mod\ \goth D^\ell_\alpha$ (as
it is $\aleph_1$-saturated), so we have many candidates for
$B^\ell_\alpha$.

Now we present the construction itself. Assume $B^0_i$, $B^1_0$ have
been choosen for $i<\alpha$ and we shall define $B^0_\alpha$, $B^1_\alpha$.
As $\langle B_i^0,B^1_i:i<\alpha\rangle$ is
defined, so is $Q_\alpha$. 

Let
$$
\align
K_\alpha\!=\!\bigl\{(\zeta,m,n,i,f',g'&,\eta)\!:\,\zeta\!<\!\alpha\, \wedge\,
m\!\leq\! t_\zeta\, \wedge\, n\!<\!\omega \, \wedge\,
i\!<\!\delta_\zeta\, \wedge\, (f',g')\!\in\! Q_\alpha\\
& \wedge\, |Dom(f')|=m\ \wedge\ \eta\in {}^m\alpha\, \wedge\,
\eta\text{ strictly increasing}\\ 
& \wedge\ Dom(f')= Dom(g')=\{\eta(0),\dots\eta(m-1)\}\\
& \wedge\ f'(\eta(s))=k^\zeta_s (\text{ for }s<m)\\
& \wedge\ g'(\eta(s))=j^\zeta_s (\text{ for }s<m)\bigr\}.
\endalign
$$

Clearly, $K_\alpha$ is countable.  Let
$L_\alpha=\{(Y^0,Y^1): Y^\ell\subseteq A^\ell_\alpha$ for
$\ell=0,1\}$, and for $(Y^0,Y^1)\in L_\alpha$, let
$Q_\alpha [Y^0,Y^1]$ be defined as $Q_{\alpha+1}$, had we
chosen $(B^0_\alpha,B^1_\alpha)$ to be $(A^0_\alpha\setminus Y^0,
A^1_\alpha\setminus Y^1)$. We say that $(Y^0,Y^1)
\in L_\alpha$ satisfies $(\zeta,m,n, i, f',g',\eta) \in
K_\alpha$ if:

\underbar{either} 

\smallskip

{\bf $(\alpha)$} for some
$\beta<\delta_\zeta$, $(f^\zeta_\beta,g^\zeta_\beta)$ and
$(f'\cup\{< \alpha, k^\zeta_m>\}, g' \cup\{<\alpha,
j^\zeta_m>\}, h)$ are compatible conditions in $Q_\alpha
[Y^0,Y^1]$ and $\beta e_{\zeta, n}i$ where $j e_{\zeta, n}i$ means:
$$
\align
(\forall \ell\leq t^\zeta)[& A^0_{\eta^\zeta_j(\ell)}\cap n =
A^0_{\eta^\zeta_i(\ell)}\cap n\ \wedge\ A^1_{\eta^\zeta_j(\ell)}\cap n =
A^1_{\eta^\zeta_i(\ell)}\cap n\\
& \qquad\wedge B^0_{\eta^\zeta_j(\ell)} \cap n =
B^0_{\eta^\zeta_i(\ell)}\cap n\ 
\wedge \ B^1_{\eta^\zeta_j(\ell)}\cap n = B^1_{\eta^\zeta_i(\ell)}].
\endalign
$$

\underbar{or} 

\smallskip

{\bf $(\beta)$} for some natural
number $u<\omega$, for every $(Z^0,Z^1) \in L_\alpha$
satisfying $Z^0 \cap u=Y^0\cap u\,\&\,Z^1\cap u=Y^1\cap
u$, clause $(\alpha)$ fails even for
$Q_\alpha[Z^0,Z^1]$.

Now, $L_\alpha$ is a complete separable metric space $\bigl($by
the metric  $d$ defined as $d((Y^0,Y^1),\,(Z^0,Z^1))=
\min((Y^0\triangle Z^0) \cup (Y^1\triangle Z^1))$ where $\triangle$ is
the symmetric difference i.e. $Y\triangle Z= (Y\setminus Z)\cup
(Z\setminus Y)$.)

Clearly:

\line{$(*)$\hfil \vtop{\hsize=0.95\hsize \noindent
 for each $(\zeta,m,n,i, f', g',\eta)\in K_\alpha$ the set
$$
L^\alpha_{(\zeta,m,n,i,f',g',\eta)}=\bigl\{(Y^0,Y^1)\in L_\alpha:
(Y^0,Y^1)
$$ 
satisfies $(\zeta,m,n,i,f',g',\eta)\bigr\}$ is an open
dense set.
}}

As $K_\alpha$ is countable, we can find a
$\langle(Y^0_\xi,Y^1_\xi):\xi<2^{\aleph_0}\rangle$ such
that:

(a) $(Y^0_\xi,Y^1_\xi) \in \cap
\{L^\alpha_{(\zeta,m,n, i,f',g',\eta)}
:(\zeta,m,n, i,f',g',\eta)\in K_\alpha\}$

(b) for $\zeta<\xi$ the set $Y^0_\zeta \cap Y^0_\xi$ is
finite.

(this is like building a perfect set of Cohen generic
reals.)

So (by the $\aleph_1$-saturation), for some $\xi< 2^{\aleph_0}$ we have
$Y^0_\xi=\emptyset\, mod\,\goth D^0_\alpha$ and $Y^1_\xi=
\emptyset\,mod\, \goth D^1_\alpha$. 

We let $(B^0_\alpha, B^1_\alpha)=(A^0\setminus Y^0,
A^1_\alpha\setminus Y^1)$.

So we have finished the inductive definition of the
$(B^0_\alpha, B^1_\alpha)$'s.  

Now we show that the construction quarantees that $Q$ satisfies the
c.c.c. Suppose $\{(f_i, g_i): i<\omega_1\}$ is an uncountable
antichain. As we explained above, we may assume there are $t_\zeta$,
$\delta_\zeta, \ldots$, and $\{(f_i, g_i): i<\delta_\zeta\}$ is an
``elementary subsystem''. So in particular, if for some $u<\omega$ and
$F^\ell: t_\zeta+1 \rightarrow {\Cal P}(u)\times {\Cal P}(u)$
(where $\ell\in \{0, 1\}$)
we have an $i<\omega_1$ such that
$(A^\ell_i\cap u, B^\ell_{\eta^\zeta_i(m)}\cap u) =
F^\ell(m)$ for each $\ell<2$ and $m\leq t_\zeta$, then such $i$ exists
already below $\delta_\zeta$.
So for each $u <\omega$ and $\gamma<\omega_1$ for some
$i_{\gamma, u}<\delta_\zeta$ we have $\gamma e_{\zeta, u}i_{\gamma, u}$.
But now consider $\gamma>\delta_\zeta$
such that $\eta^\zeta_\gamma(0)> \delta_\zeta$ and
$\eta^\zeta_\gamma(0)>\zeta$, and let $m\leq t_\zeta$ be
the least ordinal $n^\zeta_\gamma(m)$ for which condition $(\alpha)$
fails for some finite $u<\omega$ and $i=i_{\gamma, u}$ which is
$<\delta_\zeta$
(note that if $\alpha=\eta^\zeta_\gamma(t_\zeta)$ and $u=0$ then
condition $(\alpha)$ cannot possibly hold). Let it fail for $u=u_m$.
Then condition $(\beta)$
holds. But if
$u$ witnesses $(\beta)$ then for $\max\{u_m, u\}$ we have a failure for
$m-1$, contradiction.

The construction guarantees
that $Q$ satisfies the c.c.c, hence we have finished.
\null\hfill$\square_{2.2}$

\null\hfill$\square_{2.1}$
\bigskip
\subhead{2.3 Remark}\endsubhead 1) We can act as in [Wi:C 2], and then use
really $\aleph_1$-saturation and all filters are in $V$.

2) Can we ignore $CH$ and make $2^{\aleph_0}$ larger?  

Assume $\lambda=\lambda^{\aleph_1}= cf(\lambda)$ and 
$\diamondsuit_S$ where
$S=\{\delta <\lambda:cf(\delta) =\aleph_1\}$
and
$(\forall \alpha<\lambda)[|\alpha|^{\aleph_0}<\lambda]$.
We can find a forcing notion $P$, which is c.c.c. of cardinality
$\lambda$, and $\Vdash_P$ ``$2^{\aleph_0}=\lambda$ and $MA$ and $(*)$
of 2.1 holds.''

Why? we
use finite support iteration, if $\alpha \not\in
S, Q_\alpha$ is adding a Cohen real; if $\delta\in S$ and
$\diamondsuit_S$ guesses $P,\langle A^\ell_i: i
<\omega_1, \ell<2\rangle$, and $\langle \goth
D^\ell_\alpha \cap \Cal P(\omega)^{V^{P_\alpha}}:
\alpha<\omega_1,\ell< 2\rangle$ as in the proof, we imitate
the proof, but for $\langle(B^0_\alpha,
B^1_\alpha):\alpha<\omega_1\rangle$ we use a sequence
$\langle T_\alpha:\alpha<\omega_1\rangle,$ $T_\alpha$ a
perfect set of members of $L_\alpha$ such that
for every
large enough $\alpha<\omega_1$, all branches of
$T_\alpha$ are Cohen generic over $V^{P_\delta}$.
\bigskip

To answer the question of the referee:
\proclaim{2.4 Claim} The statement $\otimes$ of 2.1 follows from $PFA$.
\endproclaim

\demo{PROOF} 
Consider the forcing $P= Levy(\aleph_1, \aleph_2)* \name{Q}$
where $\name{Q}$ is constructed as in the proof of 2.1 in the universe
$V^{Levy(\aleph_1, \aleph_2)}$ (note that forcing with
$Levy(\aleph_1, \aleph_2)$ adds no
reals, so $\goth D_i$ is still an ultrafilter and $CH$ holds) and let 
$\name{B}$ be the name of the desired set. So $\Vdash_P$ ``for
$i<\omega_1$, $\ell<2$, for some $A\in \goth D^\ell_i$ we have
$[\ell=0 \Rightarrow \name{B} {}^*\supseteq A]$ and $[\ell=1
\Rightarrow \name{B}\cap A =^* \emptyset]$". Apply $PFA$.
\enddemo
\bigskip

\head{\S 3. Concerning Arhangelskii's Problem}\endhead

\proclaim{3.1 Theorem}  The following is  consistent with
$ZFC\  +\  GCH:$ 

\line{$(*)$\hfil \vtop{\hsize=0.95\hsize \noindent
  There is a regular space of cardinality $\aleph_2$
which is Lindel\"of and has pseudo-character $\aleph_0$.}}
\endproclaim
\bigskip
\subhead{Remark}\endsubhead  We had first said ``a Hausdorff space...''
but Kunen noted the proof actually yields a regular space.
\bigskip
\demo{PROOF:}  Assume $V$ satisfies $GCH$.  Let $P$ be
the set of tuples of the form
$p=\langle A,f, E, T\rangle=\langle A^p, f^p, E^p, T^p\rangle$, where:
\roster
\item $A$ is a countable subset of $\omega_2$,
\item $f$ is a two place function from $A$ to $\omega+1$,
where we write $f_x(y)$ instead of $f(x,y)$, and $f_x(y)=\omega
\Leftrightarrow y=x$ holds.
\item $E$ is a three place relation on $A$, but we write
$x E_\gamma y$, and demand that for each fixed $\gamma$, $E_\gamma$
is an equivalence relation, and $\gamma \in A\wedge
\gamma<\beta \wedge x E_\beta y \Rightarrow x E_\gamma y$,
while $x<\gamma\wedge x E_\gamma y\Rightarrow x=y$. We stipulate
$E_{\omega_2}$ as the equality on $A$.
\item $T$ is a countable set. Each member $B$ of $T$ will
be called a formal cover, $|B|=\aleph_0$, and $B$ is of the
form $\{\tau^B_n:n<\omega\}$, where each  $\tau^B_n$ is
the formal intersection of finitely many $U^n_x(x\in A,
n\in \Bbb Z)$ $\bigl[$the intended meaning is: $U^n_x=\{y:f_x(y) \geq
n\}$ for $n\geq 0$, and $U^{-n}_x = \omega_2\setminus U^n_x$ for
$n>0$; we say $p\vDash$ ``$y \in \tau$" in the natural
case  (i.e. $y\in A$, 
$\tau=\bigcap_{i<k}U^{\ell(i)}_{x(i)},$
where $k<\omega, \ell(i)\in
{\Bbb Z}, x(i) \in A$ and for each $i<k$ we have:   
$f_{x(i)}(y)\geq|\ell(i)|$ iff $\ell(i)\geq 0)\bigr]$.  We let,
 for a formal term $\tau$, $dom (\tau)=\{x\in A:x$ is
mentioned 
in $\tau\}, Dom(B)=\bigcup\{dom (\tau):\tau \in
B\}$.

The real restrictions are

\item{} if $z E_\gamma y$, and $ x<\gamma$ then
$f_x(y)=f_x(z)$.
\item\ \ \ \, (a)  if $B\in T, x\in A$, then $p\vDash``x\in
\tau^B_n$'' for some $\tau^B_n \in B$.

\indent (b) moreover, for each finite $A^1\subseteq A$ and
$h:A^1\rightarrow\omega$, there is a $\tau \in
B$ such that 
$$[x\in A^1\,\wedge\, U^n_x \text{ appears in }\tau
]\Rightarrow\left[[0\leq n\leq h(x)]\vee [n<-h(x)]\right]$$

\indent (c) moreover, for each $B\in T, x\in A$ and
$\gamma\in \{\omega_2\}\cup Dom (B)$ satisfying $\gamma \leq
x$, and a finite $A^1\subseteq A$ and function
$h:A^1\rightarrow \omega$, there is a 
$\tau=\bigcap_{i<k}
U^{\ell(i)}_{x(i)} \in B$ such that: $[x(i)<
\gamma]\Rightarrow[x\in U^{\ell(i)}_{x(i)}]$ and
 $[x(i)\geq \gamma \wedge x(i) \in A^1] \Rightarrow
\left[[0\leq \ell(i) \leq h(y)] \vee [\ell(i)
< -h(y)]\right]$.  Note that for
$\gamma=\min\ Dom (B)$, clause (c) reduces to clause (b) and for
$\gamma=\omega_2$, clause (c) reduces to clause (a). \endroster
 \enddemo 
\bigskip

\noindent
[\underbar{Explanation:} The set $A$ approximates the set of points, the
function $f$ describes the $U^n_x$'s which will generate the topology
as clopen sets, $T$ is a set of ``countable covers'', i.e. we think
of a possible covering which is a counterexample to Lindel\"ofness and
``promise'' that a countable subfamily of such cover, will cover the
entire space.
In demand (6), clause (a) just says that each $B\in T$ really covers,
clause (b) is necessary when we prove e.g. density of $\{p: x\in A^p\}$
for $x\in \omega_2$, (see the proof of Fact B). This is done by an
increasing $\omega$-sequence of descriptions of the important new
values of $f$, so this clause tells us a finite information, so does not
prevent us from preserving ``$B$ is a cover''. 

Still, why do we need clause (c) of demand (6)? We want that our
forcing notion satisfies the $\aleph_2$-c.c., so we use the
$\Delta$-system lemma, and during the construction (i.e. the proof of
Fact E i.e. the construction of a common upper bound of $p_1$, $p_2$) we have
finitely many commitements on new values of $f$, we want to make $x\in
\bigcup\limits_{n<\omega} \tau^B_n$ for $B\in T^{p_j}$, $x\in
A^{p_i}\setminus A^{p_j}$, so $f_y(x)$ is determined for all $y\in
A^{p_i}\cap A^{p_j}$ and for finitely many $y\in A^{p_j}\setminus
A^{p_i}$, and clause (6)(c) guarantees we can deal with this. We have above
avoided ``justifying'' the use of the equivalence relation $E^p$, it is
needed when in the proof of Fact D (Lindel\"ofness holds), to the union
of a generic enough sequence $\langle p_n: n<\omega\rangle$ for the
union $q$ we add to the $T^q$ a covering thus defeating a possible
counterexample, we need $E$ to verify condition (6)(c).]
\bigskip
If there are several $p$'s in consideration, we shall write $A^p,f^p,\dots$
or $p^\ell=\langle A^\ell,f^\ell,E^\ell,T^\ell\rangle$. 
Now we define the order on $P:p\leq q$ iff $A^p\subseteq
A^q, f^p=f^q \restriction A^p, E^p=E^q\restriction A^p,
T^p\subseteq T^q$.  In $V^P$   we define the following
topology on $\omega_2$. For $x\in\omega_2, n<\omega$, let
$U^n_x=\{y:f^p_x(y)\geq n$ for some $p$ in the generic set
$\}$, and if
 $n> 0$, we let $U^{-n}_x=\omega_2  \setminus U^n_x$.
Now, $\{U^n_x: x\in\omega_2, n\in {\Bbb Z}\}$ will be closed and
open, and the topology $X$ is the minimal one which
satisfies this.  So the set of finite intersections of 
$U^n_x$'s
forms a basis.  By clause (2) in the definition of $P$,
and the Fact B below, we know $\bigcap_{n<\omega}
U^n_x=\{x\}$, so as each $U^n_x$ is clopen, the space is
Hausdorff and even regular, and has pseudo-character
$\aleph_0$. \bigskip
\subhead{Fact A}\endsubhead  $P$ is $\aleph_1$-complete; in fact, any
ascending $\omega$-sequence has a naturally defined union.

In fact, we already use

\subhead{Fact B}\endsubhead  For every $ p\in P$ and $z\in \omega_2$
there is a $q\in P$ such that $q \geq p$ and $z\in A^q$.

Moreover, if $z\not\in A^p$, for any finite subset $A^*$
of $A^p$ and function $h^*$ from $A^*$ to $\Bbb Z$ we can
demand $f_x^q(z)=h^*(x)$ for $x\in A^*$.
\bigskip
\demo{PROOF OF B:} The non-trivial part is to satisfy
clause (6). We first
define $f_x(z)$ for $x\in A^p$ to satisfy (6)(c) when $z$ here stands
for $x$ there. 
So let $\{(B_k,A_k, h_k, \gamma_k):k<\omega\}$ be a list of all
tuples $(B, A, h, \gamma)$ such that $\gamma\in
\{\omega_2\}\cup Dom (B)$, $\gamma\leq z$, $B\in T^p$, $A\subseteq A^p$ is
finite and $h: A\rightarrow \omega$.  
Now we define by induction on $k$ a finite set
$D_k\subseteq A^p$ and $f_x(z)$ for $x\in D_k$. 

For
$k=0$, $D_0=A^*$, $\bigwedge_{x\in D_0} f_x(z)=h^*(x)$. 
If we have defined $D_k$,
let us define $h'_k$, $A'_k$ as follows: $A'_k= A_k \cup D_k$ and
$h'_k(x)$ is $h_k(x)$ if $x\in A_k$ and $f_x(z)$ if $x\in D_k\setminus
A_k$;
choose $\tau \in B_k$ as
exists by (6)(c) (with $A'_k$ in place of $A^1$ and $h'_k$ in place of
$h$), let $D_{k+1}= D_k \cup Dom(\tau)$, and define $f_x(z)$
for $k \in Dom(\tau)\setminus D_k$ as $\ell $ if
$U^\ell_z$ appears in $\tau$ with $\ell\geq 0$ and as $0$ otherwise.

We can at last complete the definition of $h_x(z)$ for
$x\in A \setminus \bigcup_{k<\omega}D_k$.  
Lastly define $f_z(y)$ for  $y\in A^p\cup \{z\}$ as $\omega$ if $y=z$
and 1 if $y\neq z$.
If we let
$q=\langle A\cup\{z\}$, $f^q$ is $f$ expanded as described above, $ E^p$
(i.e. $x$ is $E_\gamma$-equivalent only to itself),
$T^p\rangle$ then $q$ is O.K. \hfill$\square_{B}$

\medskip

Similarly,
we can prove

\bigskip

\subhead{Fact C}\endsubhead  
\roster
\item For every $p\in P$, $z\in \omega_2
\setminus A^p$ and $\gamma \leq y$ in $A^p,\gamma \leq z$,
there is a $q\in P$ such that $p\leq q$ and $q\vDash
``z E_\gamma y"$.  Moreover, for a given finite $A^1\subseteq
A^p \setminus\gamma$, and a function $h^1:A^1\rightarrow
\omega$, we can demand $q\vDash ``f_x(z)=h^1(x)$ for $x\in
A"$. 
\item The following is a dense subset of $P$, closed under unions:
$$
\align
{\Cal I}=\biggl\{q\in P: & \text{ \underbar{for every} } \gamma\in
A^q\cup\{\omega_2\}\text{ and finite }A\subseteq A^q\text{ and}\\
& \text{ function }h: A\rightarrow \omega\text{ and }x\in A^q\text{
satisfying }\gamma\leq x\\
&\text{ \underbar{there is} }x'\in A^q\text{ such that
}\gamma\leq x', x'E_\gamma x \text{ and }\\
&\qquad (\forall y\in A\setminus \gamma)[f_x(y)=h(y)]\biggr\}.
\endalign
$$ 
\endroster
\bigskip
\demo{PROOF OF C:} 1) like the proof of Fact B.

2) For any $p\in P$, choose $p_n\in P$, $p_0=p$, $p_n\leq p_{n+1}$,
each time use part (1) with a suitable bookkeeping and take the union by
Fact A.
\hfill$\square_{B}$
\enddemo
\bigskip 
\subhead{Fact D}\endsubhead
  The space is Lindel\"of.
\bigskip

\demo{PROOF OF D:} Let $\name\sigma$ be a name of a cover and 
$\name\sigma_x$ a member of it  to which $x$ belongs, and
$w.l.o.g.\ \name{\sigma}_x$ is a member of the basis (i.e. is a
$\tau$).  Now, for each $p\in P$ and $x\in A$, there is a $q \geq
p$, such that $q\Vdash ``\name\sigma_x=\tau"$ for some
specific $\tau$.  By Fact A for every $p\in P$ there is $q\in P$ such
that 
$$
p\leq q\in {\Cal I}\ \wedge\ \bigwedge_{x\in A^p}\bigvee\{q\Vdash
``\name{\sigma}_x=\tau": \tau\text{ as in demand (4)}\}.
$$
So for every $p\in P$ we can find $\langle p_n: n<\omega\rangle$ such
that $p\leq p_0$, $p_n\leq p_{n+1}$, $p_n\in {\Cal I}$ and $x\in A^{p_n}
\Rightarrow  p_{n+1}\Vdash~~\name{\sigma}_x=\tau^x"$. Let
$q=\bigcup_{n<\omega} p_n$, now $q\in P$ is an upper bound of $\{p_n:
n<\omega\}$ and $q\in {\Cal I}$ by Fact C(2). Now
$$
B^*=\{\tau:q\Vdash ``\name\sigma_x=\tau\text{'' for some }x\in
A^q \text{ and }\tau\text{ as in (4) (for } q\text{)}\},
$$ 
satisfies the requirements on B in clause (4). Now $B^*$ also satisfies
the requirements of part (6) in the definition: clause (a) holds as
$$
x\in A^q \Rightarrow \bigvee_n x\in A^{p_n} \Rightarrow \bigvee_n
(p_{n+1}\text{ forces a value to }\name{\sigma}_x).
$$
Clause (c) holds by the above as $q\in {\Cal I}$.
So
$q^*\buildrel\text{df}\over =\langle A^q,f^q,E^q,
T^q\cup\{B^*\}\rangle \in P$ [why? check; the main point is (6)(c)
which holds as $q\in {\Cal I}$]. Also $q\leq q^*$, so $q^*$ forces that $B^*$
is (essentially) a countable subcover, as required. \hfill$\square_{D}$
\enddemo 
\bigskip
\subhead{Fact E}\endsubhead  $P$ satisfies the $\aleph_2$-chain
condition.  (Hence in $V^P,G$ has power $>2^{\aleph_0}$).
\bigskip

\demo{PROOF:}
Let $p_i\in P (i<\omega_2)$.  It is well known that we can
assume that for some $p$ and
$\alpha_i\,(\text{for} \  i<\omega_2)$: $\alpha_i $ is
increasing,  $A^{p_i}\cap \alpha_i=A^p, A^{p_i} \subseteq
\alpha_{i+1}, p\leq p_i$. Now, like in the proof of Fact B, we
can prove $p_0, p_1$ can be extended to a condition
$\langle A^{p_0} \cup A^{p_1}, f^{p_0} \cup f^{p_1},
E^{p_0} \cup E^{p_1}, T^{p_0} \cup T^{p_1}\rangle$.
\hfill$\square_{E}$\newline
\null\hfill$\square_{3.1}$
\enddemo

\subhead{Remarks}\endsubhead
\roster
\item The proof works with $\aleph_0$ replaced by any
$\lambda$ such that $\lambda^{<\lambda}=\lambda,
2^\lambda=\lambda^+$; countable is replaced by ``of cardinality $\leq
\lambda$'' and $\tau$ are still finite formal intersections,
$Rang(f)\leq \omega=1$.  We get $\lambda$-Lindel\"of but
still pseudo-character $\aleph_0$.
\item  It is well known that there is no Lindel\"of space
of pseudo-character $\aleph_0$ of power $\geq$ `` first
measurable''. \endroster
\bigskip

\head{\S 4. More on Arhangelskii's problem}\endhead

We prove:

\proclaim{4.1 Theorem} If ($ZFC + GCH +\exists$ a weakly compact) is
consistent, then so is the following: $ZFC + CH +
not(*)$ where

\line{$(*)$\hfil \vtop{\hsize=0.95\hsize \noindent
 There is a regular space of cardinality $\aleph_2$
which is Lindel\"of and has pseudo-character $\aleph_0$.}}
\endproclaim
\bigskip
\subhead{Question}\endsubhead ($CH$) Is there a Lindel\"of space with 
pseudo-character $\aleph_0$ and with $2^{\aleph_1}$ points? More than
$\aleph_2$ points?
\bigskip
\demo{PROOF:}  
Let $V\vDash ``GCH\wedge \kappa$ weakly compact.''

Let $P_i$ ($i<\kappa^+$) be the forcing for adding a Cohen subset to
$\omega_1$ (so $P_i$ is $\aleph_1$-complete).

$Q^\kappa_2$ is the Levy collapse of $\kappa$ to $\aleph_2$
(i.e. every $\theta \in (\aleph_1,\kappa)$ is collapsed to
$\aleph_1$, and each condition is countable).

$Q_1=\Pi_{i<\kappa^+} P_i=\{p\in \Pi_{i<\kappa^+}P_i:
p$ has countable support (i.e. $p(i)=\emptyset$ for all except
countably many $ i)\}$.

$Q=Q^\kappa_2 \times Q_1$.

In $V^Q$ we know $\kappa$ is $\aleph_2$, $CH$ holds,
$2^{\aleph_1}=\kappa^+$ (and only the cardinals
$\theta\in(\aleph_1,\kappa)$ were collapsed.)

We prove that in $V^Q$ there is no Lindel\"of space $X$ with
countable pseudo-character, such that $|X|=\kappa$.

If there is such an $X$, we can assume its set of points is
$\omega_2$, and $x\in U^n_x, U^n_x$ open,
$\bigcap_n U^n_x=\{x\}$ for all $x\in\omega_2$
(by the countable pseudo-character).

So the topology is the minimal one to which all $U^n_x$ belong; this
is O.K. as if we decrease family of open sets, the Lindel\"ofness is
preserved. 

We $w.l.o.g.$ identify $X$ with $\langle(x, n, U^n_x): x<\omega_2,
n<\omega\rangle$. 
So the topology is the minimal one to which all $U^n_x$ belong; this
is O.K. as if we decrease the family of open sets, the Lindel\"ofness is
preserved. Note

\item{$\oplus$} if $X\in V_1\subseteq V_2$ and $U\subseteq X$ with
$U\in V_1$, then
$V_1\vDash$ ``$U$ is open in $X$'' iff $V_2\vDash$ ``$U$ is open in
$X$''.

(because $U$ is open iff $\forall x\in U 
\exists x_1, \ldots, x_k \exists n_1,
\ldots, n_k (x\in \bigcap\limits_{\ell=1}^k U^{n_\ell}_{k_\ell} \subseteq U$))
\enddemo
\bigskip
\proclaim{Claim A}  In $V^Q$ (but $CH$ is used), for some closed unbounded
$C\subseteq \kappa$, we have  $\alpha \in C \wedge cf(\alpha)
>\aleph_0\Rightarrow (X\restriction \alpha$ is not
Lindel\"of, moreover, there is a 
$g:\alpha\rightarrow\omega$, so that $\{U^{g(x)}_x:x\in
\alpha\}$ has no countable subcover).
\endproclaim
\bigskip
\demo{PROOF OF A}:  $C$ will be the family of
$\alpha<\omega_2$ such that:

{\it if} $Dom (h)$ is a countable bounded subset of
$\alpha$,  $Range (h)\subseteq \omega,
\omega_2\not=\bigcup\{U^{h(x)}_x: x\in Dom (h)\}$
{\it then\/} there is a $\beta<\alpha,\, \beta
\not\in\bigcup\{U^{h(x)}_x:x\in Dom (h)\}$.

Clearly, $C$ is
closed and by $CH$ it is unbounded. 

If $cf(\alpha) >\aleph_0$, we can omit the ``bounded'', as
every countable subset of $\alpha$ is bounded.

For $\alpha\in C$ such that $cf(\alpha)>\aleph_0$, define
$g:\alpha\rightarrow\omega$ as follows $: g(x)$ is the
first $n<\omega$ such that $\alpha\not\in U^n_x$ (exists as
$\bigcap_{n<\omega}U^n_x=\{x\})$.

Clearly $\{U^{g(x)}_x: x<\alpha\}$ cover $\alpha$.
Suppose $\{U^{g(x)}_x:x < \alpha\}$ has a countable
subcover $\{U^{g(x)}_x:x\in Y\}, |Y| \leq \aleph_0$.  Let
$h=g \restriction Y$ and we get a contradiction to the
definition of $C$ (because $\alpha$ witnesses that the
union is not $\omega_2$).
\hfill$\square_A$
\enddemo
 \bigskip 
\proclaim{Claim B} Suppose $V$
satisfies $CH$ and $X$ and $U^n_x(x\in X)$ are as above.  Suppose
$X$ is Lindel\"of, $P$ is an $\aleph_1$-complete forcing, but
in $V^P$, the space $X$ is not Lindel\"of.

Then also in $V^{P_0}$ (remember that $P_0$ is adding one
Cohen subset to $\omega_1)$, the space is not Lindel\"of.
\endproclaim
\bigskip
\demo{PROOF OF B:}  Suppose $\name{\tau}$ is a $P$-name for a
cover contradicting Lindel\"ofness, $wlog$ the cover
consists of old sets.  Let $p\in P$.

We define $p_\eta\in P$ for $\eta
\in{}^{{\omega_1}>}\omega$ by induction on the length $\ell(\eta)$
of $\eta$,
an old open set $U_\eta$ (where $\ell(\eta)$ is a successor) such that:
\roster
\item $p_{\eta\restriction\alpha} \leq p_\eta$ ($\le$ means ``weaker than")
\item $p_\eta\Vdash_P ``U_\eta\in\name{\tau}"$ (for $\ell(\eta)$ a
successor)
\item $X=\bigcup_{n<\omega}U_{\eta\hat{\,}\,\langle
n\rangle}$.
\endroster
\bigskip
For $\ell(\eta)=0, \eta=\langle \  \rangle, p_\eta=p$, for
limit only (1) applies and we use $\aleph_1$-completeness.  If 
$\eta\in (^\alpha \omega)$ and 
$  p_{\eta\restriction(\beta+1)}$ is defined for
$\beta<\alpha$, let
$$
\align
F_\eta=\{U: & \ U\text{ an open set of }X\text{ (in
the universe }V\text{)}\\
& \text{\ \ \ and for some }q\geq p_\eta, q\Vdash ``U\in
\name{\tau}"\}.
\endalign
$$
 Clearly, $F_\eta$ is a cover.

$F_\eta$ is a cover, but $X$ is Lindel\"of, so for some
countable $F'\subseteq F_\eta$ we have $X=\bigcup F'$.  Let
$F'=\{U^\eta_n:n<\omega\}$ (maybe with repetitions).  Let
$U_{\eta\hat{\,}\,\langle
n\rangle}\mathop{=}\limits^{\text{def}}U^\eta_n$, 
and let
$p_{\eta\hat{\,}\,\langle n\rangle} \geq p_\eta$ be chosen
so that $p_{\eta\hat{\,}\,\langle n\rangle}\Vdash_P$
``$U_{\eta\hat{\,}\,\langle n\rangle}\in \name{\tau}$.''

Now we show that in $V^{P_0}$ the space $X$ is not
Lindel\"of.  For a generic $g\in (^{\omega_1} 2)$
let $\sigma$ be the family 
$\{U_{g\restriction(\alpha+1)}:\alpha<\omega_1\}$. It is easily seen
that $\sigma$ is a cover of $X$. Suppose $X\subseteq
\bigcup\{U_{g\restriction (\alpha+1)}: \alpha<\beta\}$ for some
$\beta<\omega_1$. Without loss of generality, $\beta$ is a limit
ordinal. Then $p_{g\restriction \beta}\Vdash_P ``\sigma\subseteq
\name{\tau}$'', in contradiction with the choice of $\name{\tau}$.
\hfill$\square_B$
\enddemo

\demo{Continuation of the proof of 4.1:}
W.l.o.g. $y\in U^n_x$ is determined in $Q_3= Q^\kappa_2
\times \Pi_{i<\kappa} P_i$. [Why?
As $Q$ satisfies the $\kappa^+$-c.c. there are $\kappa$ maximal
antichains ${\Cal I}_{y, n, x}$ of elements forcing a truth value to
``$y\in U^n_x$''. So $|\bigcup_{y,n,x} {\Cal I}_{y, n, x}|\leq \kappa$, so
for some $\alpha$ we have $\bigcup{\Cal I}_{y, n,x}\subseteq Q^\alpha_2 \times
\prod_{i<\kappa}P_i$.] Also 
in $V^{Q_3}$, the space $X$ (i.e. the space defined by letting 
a subset be open iff
it is forced to be open, see $\oplus$ above) is Lindel\"of, 
noting that a cover in $V^{Q_3}$ is also a cover in $V^Q$ and those
two universes have the same $\omega$-sequences of members of $V^{Q_3}$
as no $\omega$-sequences
are added by $\Pi_{i\in[\kappa,\kappa^+)}P_i$).  So,
forcing by $P_0=P_\kappa$ over $V^{Q_3}$ does not contradict
Lindel\"ofness, by Claim B.

As $\kappa$ is weakly compact, for some  stationary set $S$ we have
$$S\subseteq\{\alpha<\kappa: \alpha\text{ strongly inaccessible}\},$$ 
and in addition, for each $\alpha\in S$, 
we can split the forcing $Q$ to 
$Q'_\alpha \times Q''_\alpha$, both $Q'_\alpha$ and $Q''_\alpha$ are 
$\aleph_1$-complete forcings,
and
$Q'_\alpha=Q_2^\alpha \times\Pi_{i <\alpha} P_i$ 
so that:
$$Q'_\alpha\text{ determines} \ ``y\in U^n_x" \text{ for }\ 
y,x<\alpha, n<\omega,$$
and in $V^{Q'_\alpha}, \alpha$ becomes $\aleph_2$. Also,
the part of the space $X$ that we get after forcing
with $Q_\alpha^{'}$, $X_\alpha$ that is $X\restriction \alpha$, is
Lindel\"of of pseudo-character $\aleph_0$, as exemplified by the $U^n_x$, and
adding a $P_0$-generic does not contradict Lindel\"ofness. 
(Here we use the weak compactness of $\kappa$ i.e.
$\Pi^1_1$-indescribability of $\kappa$.)

Now, by claim A, for some such $\alpha$, in $V^{{Q'_\alpha}
\times {Q''_\alpha}}$ the space $X_\alpha$ is no longer
Lindel\"of. Therefore, forcing by $Q''_\alpha$ abolishes
Lindel\"ofness. Also $Q'_\alpha$, $Q''_\alpha$ are
$\aleph_1$-complete, so $Q''_\alpha$ is $\aleph_1$-complete in
$V^{G'_\alpha}$, hence (Claim B) $P_0$ forcing abolishes Lindel\"ofness,
a contradiction.
\hfill$\square_{4.1}$
\enddemo
\bigskip

\subhead{4.2 Remark} \endsubhead
Note that during the proof we did not use the regularity of $X$.

\bigskip

{\obeylines
\leftline{\vtop{\hsize=0.40\hsize
Institute of Mathematics 
The Hebrew University 
Jerusalem, Israel
}\ \ \vtop{\hsize=0.40  \hsize
Rutgers University 
Department of Mathematics 
New Brunswick, NJ USA}}}

\vfill
\pagebreak

\shlhetal
\end
\Refs
\ref \by [B] J.Baumgartner, {\bf All $\aleph_1$-dense sets
of reals can be isomorphic}, {\it Fund. Math.\/},
1973 vol. 79 (101--106)\endref

\ref \by [Go] I. Gorelic, {\bf The Baire category and forcing large
Lindel\"of spaces with points $G_\delta$}, {\it Proceedings of Amer.
Math. Soc.}, 1993 vol. 118 \endref 

\ref \by [HJ] A. Hajnal and I.Juhasz, 
{\bf  Lindel{\"o}f spaces \`a la Shelah}, {\it Coll\. Math\. Soc\. J\.
Bolyai, Topology}, 1978 vol. 23\endref   

\ref \by [ShSt:167] S. Shelah and L. Stanley, {\bf $S$-forcing. IIa.
Adding diamonds and more applications: coding sets, Arhangelskii's
problem and ${\Cal L}[Q^ {<\omega}_ 1,Q^ 1_ 2]$}, {\it Israel Journal
of Mathematics}, 1986 vol. 56 (1--65)\endref

\ref \by [Wi:C2 ]  E. Wimmer, {\bf The Shelah $P$-point
independence theorem}, {\it Israel Journal of Mathematics,\/}
1982 vol. 43 (28--48)\endref


\endRefs
